\documentclass[10pt]{amsart}
\usepackage[utf8]{inputenc}
\usepackage[T2A]{fontenc}
\usepackage{amsfonts,amssymb,mathrsfs,amscd,comment,bm}

\newtheorem{thm}{\sc Тheorem}
\newtheorem{prop}{\sc Proposition}

\def \Lip {\mathop{\rm Lip}\nolimits}

\begin{document}

\title{{
On bi-Lipschitz isomorphisms of self-similar Jordan arcs
}}

\author{Ilya Galai, Andrei Tetenov}
\date{}

\begin{abstract}
 We find  the conditions for bi-Lipschitz equivalence of self-similar Jordan arcs which are the attractors of self-similar zippers.
 \end{abstract}

\maketitle

\makeatletter
\renewcommand{\@makefnmark}{}
\makeatother
\footnotetext{This work was supported by the Mathematical Center in Akademgorodok with the Ministry of Science and Higher Education of the Russian Federation, agreement no. 075-15-2019-1613.}


 A system $\mathcal S =\{S_i, i\in J\}$  (where $J=\{1,...,m\}$)  of contracting similarities of $\mathbb R^n$, is
 called a {\em self-similar zipper} with vertices $\{z_0, \ldots, z_m\}$ and signature ${\bm\varepsilon} \in \{0, 1\}^m$, if for any $i\in J$, $S_i
(z_0) = z_{i-1+\varepsilon_i}$ and $S_i (z_m) = z_{i-{\varepsilon_i}}$ \cite{ATK}.
 A compact set $K\subset \mathbb R^n$ is called the {\em  attractor of the zipper} $\mathcal S$, 
if $K=S_1(K)\cup....\cup S_m(K)$.

If $0=t_0 < t_1 <...< t_m=1$  is  a set of points of the segment $I=[0,1] \subset
\mathbb R$, 
then a self-similar zipper ${\mathcal T} = \{ T_1,..., T_m\} $ with vertices $\{t_0,...,t_m\}$ and 
signature  ${\bm\varepsilon} \in \{0, 1\}^m$ is called {\em linear}.
The attractor of a linear zipper is the segment $I$. \\

For any zipper ${\mathcal S} = \{ S_1,..., S_m\} $
with vertices $\{z_0, \ldots, z_m\}$ and for any linear zipper with the same signature $\bm \varepsilon$ there is a unique continuous map $g : I \to K({\mathcal
S})$, 
for which $g(t_i) = z_i$ and $S_i\circ g =
g\circ T_i$ for every $i\in J$. 
Moreover, the mapping $g$ is   H\"older continuous  and $g(I) = K({\mathcal S})$.
Such mappings $g$ are called
 {\it structural parametrizations} of the attractor of the zipper  ${\mathcal S}$.
A zipper ${\mathcal S}$ is called a {\it Jordan zipper} if and only if
 one (and therefore any) of the structural parametrizations  of its
 attractor is a homeomorphism of the segment
$ I = [0,1] $ to $K({\mathcal S})$\cite{ATK}.
Indeed, as it was shown by Z.-Y. Wen and L.-F. Xi \cite{WX}, it is possible that the homeomorphism $ g^{- 1}$, inverse to the structural parametrization, does  not satisfy the H\"older condition for any value of the H\"older exponent.\\

Note an important property of the subarcs of the attractor $\gamma$ of a Jordan zipper $\mathcal S$. Let $J^*=\{{\bm j}=j_1...j_n, j_k\in J, n\in \mathbb N\}$ be the set of all multi-indices over $J$.
Let $S_{\bm j}=S_{j_1}\circ....\circ S_{j_n}$, аnd 
$\gamma_{\bm j}=S_{\bm j}(\gamma)$. 
The systems $\{\gamma_{\bm j},{\bm j}\in J^n,n\in \mathbb N\}$ 
form a refining sequence of partitions of the arc $\gamma$, in 
which $\gamma_{\bm i}\supseteq  \gamma_{\bm j}$ 
if and only if ${\bm i}\sqsubseteq{\bm j}$.
If ${\bm i}\not\sqsubseteq{\bm j}$ and ${\bm i}\not\sqsupseteq{\bm j}$, then $\gamma_{\bm i}\cap  \gamma_{\bm j}$ is either empty or is the common endpoint of subarcs $\gamma_{\bm i}$ and   $\gamma_{\bm j}$.\\

 We say that a Jordan  arc $\gamma\subset \mathbb R^n$ is  an {\em arc of bounded turning } \cite[p.100]{TV}, if there exists $ M>0$ such that for any $x,y\in \gamma$ the diameter $|\gamma_{xy}|$ of the subarc $\gamma_{xy}\subset\gamma$ 
with endpoints $x,y$  
is not greater than $M\|x-y\|.$ \\

Let  $\mathcal S =\{S_i, i\in J\}$ 
and $\mathcal S'=\{S'_i, i\in J\}$ be self-similar Jordan zippers, and $\gamma$ 
and $\gamma'$ be their attractors.
A homeomorphism $f:\gamma\to \gamma'$ {\em agrees} with the systems  $\mathcal S$ and $\mathcal S'$, if for any $i\in J$, $f\circ S_i=S'_i\circ f$.
In this case, we say $f$  defines an  isomorphism of the zippers $\mathcal S$ and $\mathcal S'$. Note that for any ${\bm j}\in J^*$, $f(\gamma_{\bm j})=\gamma'_{\bm j}$. \\
The uniqueness and inversibility of the structural parametrization of  Jordan zippers imply the following proposition.\\
\begin{prop} Let  
self-similar Jordan zippers $\mathcal S$ and $\mathcal S'$ with attractors $\gamma$ and $\gamma'$, 
have the same signature. Then there is a unique homeomorphism $f:\gamma\to \gamma'$, 
which agrees with $\mathcal S$ and $\mathcal S'$. $\blacksquare$
\end{prop}

We say that the zippers $\mathcal S$ and $\mathcal S'$ 
{\em are   bi-H\"older (resp. bi-Lipschitz) isomorphic} if the homeomorphism $f$ is bi-H\"older (resp. bi-Lipschitz).\\
We  prove the following statement:\\

\begin{thm} 
 Let self-similar  zippers $\mathcal S =\{S_i, i\in J\}$, $\mathcal S'=\{S'_i, i\in J\}$ have  attractors $\gamma$ and $\gamma'$ which are Jordan arcs of bounded turning. Let be $p_i$ (resp. $q_i$) be the similarity ratios  of  the mappings $S_i$(resp. $S_i'$).Then $\mathcal S$ and  $\mathcal S'$ are { bi-H\"older  } isomorphic if and only if their signatures are equal. Moreover, the Hölder exponents of homeomorphisms $f$ and $f^{-1}$ are greater or equal to $\alpha=\min\left\{{\log{p_i}}/{\log{q_i}},{\log{q_i}}/{\log{p_i}},i\in J \right\}$.
\end{thm}

\textit{Proof}. Let $p_{min} = \min{\{p_1, ..., p_m\}}$. Let
$f:\gamma\to\gamma'$ be a homeomorphism which agrees with $\mathcal S$ and $\mathcal S'.$ Let $x,y\in\gamma$, and $x'=f(x),y'=f(y)$.
There is  $M>0$ such that the arcs $\gamma_{xy}\subset\gamma$ and  $f(\gamma_{xy})=\gamma'_{x'y'}$   obey the inequalities  \begin{equation}\label{cBT}|\gamma_{xy}| \leq M\|x-y\|\mbox{  and }|\gamma'_{x'y'}| \leq M\|x'-y'\|.\end{equation}

Due to the properties of the families $\{\gamma_{\bm j},{\bm j}\in J^*\}$ and $\{\gamma'_{\bm j},{\bm j}\in J^*\}$ there are only two possibilities for the arcs $\gamma_{xy},\gamma'_{xy}$ :\\
{\bf 1.} 
There are    $\bm i=i_1...i_k\in J^*$ 
and   $i_{k+1}\in I$ such  
that  $\gamma_{{\bm i}i_{k+1}}\subset\gamma_{xy}\subset\gamma_{\bm i}$ and
$\gamma'_{{\bm i}i_{k+1}}\subset\gamma'_{x'y'}\subset\gamma'_{\bm i}$. These inclusions imply the
 inequalities $p_{{\bm i}i_{k+1}}|\gamma|\le|\gamma_{xy}|\le p_{\bm i}|\gamma|$ and
$q_{{\bm i}i_{k+1}}|\gamma'|\le|\gamma'_{x'y'}|\le q_{\bm j}|\gamma'|$. 
From the conditions \eqref{cBT} it follows that $p_{{\bm i}}\le\dfrac{|\gamma_{xy}|}{|\gamma|p_{min}}\le\dfrac{M\|x-y\|}{|\gamma|p_{\min}}$.
Since for any multi-index $\bm k$,  $q_{\bm k}\le p_{\bm k}^\alpha$,\   we get
 \[\|x'-y'\|\le|\gamma'_{x'y'}|\le p_{\bm i}^\alpha |\gamma'|\le \dfrac{M^\alpha|\gamma'|} {p_{\min}^\alpha|\gamma|^\alpha}\|x-y\|^\alpha.\]
{\bf 2.} 
There are  multi-indices $\bm i=i_1...i_k$, $\bm j=j_1...j_l$
and indices $i_{k+1},j_{l+1}$ such that
$\gamma_{{\bm i}i_{k+1}}\cup \gamma_{{\bm j}j_{l+1}}\subset\gamma_{xy}\subset\gamma_{\bm i}\cup \gamma_{\bm j}$, аnd $\gamma_{{\bm i}i_{k+1}}\cap \gamma_{{\bm j}j_{l+1}}=\gamma_{\bm i}\cap \gamma_{\bm j}$ is a singleton; similar relations hold for $\gamma'_{x'y'}$.\\ 

From the inequalities $\max(p_{{\bm i}i_{k+1}},p_{{\bm j}j_{l+1}})|\gamma|\le|\gamma_{xy}|\le(p_{\bm i}+p_{\bm j})|\gamma|$ and
$|\gamma'_{x'y'}|\le (q_{\bm i}+q_{\bm j})|\gamma'|$ 
we deduce that
$\|x'-y'\|\le|\gamma'_{x'y'}|\le 2 \max(p_{{\bm i}},p_{{\bm j}})^\alpha|\gamma'|$. As
 $\max(p_{{\bm i}},p_{{\bm j}})\le\dfrac{M\|x-y\|}{|\gamma|p_{\min}}$, we get
\begin{equation}\label{Hol}\|x'-y'\|\le 2\dfrac{M^\alpha|\gamma'|} {p_{\min}^\alpha|\gamma|^\alpha}\|x-y\|^\alpha.\end{equation}
If an index $i\in J$ 
is such that $\alpha=\log q_i/\log p_i$, then for pairs of points
$x=S_i^k(z_0), y=S_i^k(z_m)$ we have the equality $\|x'-y'\|=M_0\|x-y\|^\alpha$, where $M_0=\|z'_m-z'_0\|/(\|z_m-z_0\|)^\alpha$.
The same argument is valid for the mapping $f^{-1}:\gamma'\to\gamma$.
Thus, the mappings $f$ and $f^{-1}$ satisfy the H\"older condition with exponent $\alpha$ and this value $\alpha$ is minimal. $\blacksquare$\\

Considering the case when $\alpha=1$, we obtain the condition for  bi-Lipschitz equivalence of self-similar zippers.
 
 \begin{thm}
    Let $\mathcal{S} = \{S_1, ..., S_m\}$ and $\mathcal{S'} = \{S'_1, ..., S'_m\}$ be  self-similar Jordan zippers in $\mathbb{R}^n$  with attractors $\gamma$ and $\gamma'$, signatures $\bm \varepsilon$  and $\bm \varepsilon'$, and similarity coefficients $\Lip S_i,\Lip S_i', i\in J$ respectively.\\ 
    1. 
If  $\bm \varepsilon=\bm \varepsilon'$
and for any $i\in J$,  $\Lip S_i=\Lip S_i'$, and $\gamma$, $\gamma'$ are the arcs of bounded turning, then $\mathcal{S}$ and $\mathcal{S'}$ are bi-Lipschitz isomorphic.\\
    2. If  $\mathcal{S}$ and $\mathcal{S'}$ are bi-Lipschitz isomorphic, and $\gamma$ is an
arc of bounded turning, then $\bm \varepsilon =\bm \varepsilon',$
and for any $i\in J$,  $\Lip S_i=\Lip S_i'$, аnd $\gamma'$ is of bounded turning. 
 \end{thm}

\bigskip

{\bf  A.\,V.~Tetenov}, Novosibirsk state university. {\it E-mail}: a.tetenov@gmail.com

{\bf I.\,N.~Galay}, Novosibirsk state university. {\it E-mail}: la.a2010@yandex.ru

\end{document}